\documentclass[twoside,11pt,a4paper,leqno]{article}
\usepackage{amsfonts}
\usepackage{amsfonts}
\usepackage{amssymb,amsmath,amscd,euscript,verbatim,array}
\usepackage[center,pagestyles]{titlesec}
\usepackage{geometry}
\usepackage{anysize}
\usepackage{fancyhdr}
\usepackage{indentfirst}

\marginsize{3.5cm}{3.5cm}{2cm}{2cm}

\titleformat{\section}{\centering\normalsize}{\thesection.}{0.5em}{}
\titleformat{\subsection}{\normalsize\bfseries}{\thesubsection.}{0.5em}{}
\titleformat{\subsubsection}{\normalsize\bfseries}{\thesubsubsection.}{0.5em}{}
\newcommand{\N}{\mathbb{N}}
\newcommand{\Z}{\mathbb{Z}}
\newcommand{\R}{\mathbb{R}}
\newcommand{\Q}{\mathbb{Q}}

\newtheorem{Theorem}{Theorem}[section]

\newtheorem{Lemma}[Theorem]{Lemma}
\newtheorem{Exercise}[Theorem]{Exercise}
\newtheorem{Proposition}[Theorem]{Proposition}
\newtheorem{Remark}[Theorem]{Remark}
\newtheorem{Corollary}[Theorem]{Corollary}

\newcommand{\eps}{\varepsilon}

\newcommand{\T}{\mathbb{T}}
\newcommand{\bthm}{\begin{Theorem}}
\newcommand{\ethm}{\end{Theorem}}
\newcommand{\bpr}{\begin{Proposition}}
\newcommand{\epr}{\end{Proposition}}
\newcommand{\blm}{\begin{Lemma}}
\newcommand{\elm}{\end{Lemma}}
\newcommand{\bex}{\begin{Exercise}}
\newcommand{\eex}{\end{Exercise}}
\newcommand{\be}{\begin{equation}}
\newcommand{\ee}{\end{equation}}
\newcommand{\beal}{\begin{aligned}}
\newcommand{\enal}{\end{aligned}}
\newcommand{\brm}{\begin{Remark}}
\newcommand{\erm}{\end{Remark}}
\newcounter{item}[section]

\newcommand{\Proof}{\textbf{Proof}\hspace{0.3cm}}
\newcommand{\End}{\ensuremath{\hfill{\Box}}\\}
\renewcommand{\title}[1]{\begin{center}\textbf{\normalsize #1}\end{center}}
\renewcommand{\author}[1]{\begin{center}\small #1\end{center}}
\renewcommand{\date}[1]{\begin{center}#1\end{center}}

\makeatletter \@addtoreset{equation}{section}
\makeatother
\begin{document}
\vspace{10pt}
\title{VARIATIONAL DESTRUCTION OF INVARIANT CIRCLES}
\vspace{6pt}
\author{\sc Lin Wang}
\vspace{10pt} \thispagestyle{plain} 
\begin{quote}
\small {\sc Abstract.} We construct a sequence of generating
functions $(h_n)_{n\in\N}$, arbitrarily close to an integrable
system in the $C^r$ topology with $r<4$ for $n$ large enough. With
the variational method, we prove that for a given rotation number
$\omega$ and $n$ large enough, the exact monotone area-preserving
twist maps generated by $(h_n)_{n\in\N}$ admit no invariant circles
with rotation number $\omega$.
\end{quote}
\begin{quote}
\small Key words. invariant circles, minimal configuration,
Peierls's barrier
\end{quote}
\begin{quote}
\small AMS subject classifications (2000). 58F27, 58F05, 58F30,
58F11 \end{quote} \vspace{20pt}


\section{\sc Introduction}

 M.Herman constructed an example in [H2]. With geometrical
method he proved the example has the property that the invariant
circle (i.e. a homotopically non-trivial invariant curve) with a
given rotation number can be destructed by an arbitrarily small
perturbation in the $C^{3-\epsilon}$ topology for the exact monotone
area-preserving twist map. The KAM theory (see [S]) implies the
persistence of invariant circle with a Diophantine rotation number
under the small perturbation in the $C^{3+\epsilon}$ topology.
Hence, Herman's result is optimal. However, the perturbation in the
example in [H2] is too artificial. In this paper, following the
ideas and techniques developed by J.N.Mather in the series of papers
[M1], [M2], [M3] and [M4], we construct an example with a more
natural perturbation to achieve the same goal. More precisely, our
example has the property as follow: \vspace{6pt}

\noindent\textbf{Property:} {\it For a given rotation number
$\omega$ and $n$ large enough , there exists a sequence of
generating functions $(h_n)_{n\in\N}$, arbitrarily close to an
integrable system in the $C^r$ topology with $r<4$ such that the
exact monotone area-preserving twist maps generated by
$(h_n)_{n\in\N}$ admit no invariant circles with rotation number
$\omega$.} \vspace{6pt}

 In
[M1], [M2], [M3] and [M4], Mather introduced a notion called
Peierls's barrier as follows
\[P_\omega^h(\xi)=\min_{x_0=\xi}\sum_I(h(x_i,x_{i+1})-h(x_i^-,x_{i+1}^-)),\]
where $I=\Z$, if $\omega$  is not a rational number, $I=\{0,..., q-1
\}$, if $\omega=p/q$, and $(x_i)_{i\in I}\in \prod_{i\in
I}[x_i^-,x_i^+]$ satisfying $x_0=\xi$. Moreover, he proved that
$P_\omega^h(\xi)$ is a non-negative Lipschitz function with respect
to the variable $\xi\in\R$ with the modulus of continuity with
respect to $\omega$ and found a criterion of the existence of
invariant circles. Namely, the exact area-preserving monotone twist
map generated by $h$ admits an invariant circle if and only if
$P_{\omega}^h(\xi)\equiv 0$ for every $\xi\in \R$.

For our example, the modulus of continuity can be improved due to
the hyperbolicity of the perturbation, which follows from similar
ideas of [F]. More precisely, (Lemma \ref{apprx} below) if $\omega$
is suitable small, then
\[|P_{\omega}^{h_n}(\xi)-P_{0^+}^{h_n}(\xi)|\leq
C\exp(-n^{\delta}),\] where $\delta$ is a small positive constant
independent of $n$. Based on the improvement, we obtain that there
exists $\tilde{\xi}\in \R$ such that
$P_{\omega}^{h_n}(\tilde{\xi})\neq 0$ for $n$ large enough. It
follows that $h_n$ doesn't admit any invariant circles for $n$ large
enough.

\section{\sc Preliminaries}

\subsection{Minimal configuration}
Let $f$: $\T\times\R\rightarrow \T\times\R$ ($\T=\R/\Z$) be an exact
area-preserving monotone twist map and $h$: $\R^2\rightarrow\R^2$ be
a generating function for the lift $F$ of $f$ to $\R^2$, namely $F$
is generated by the following equations
\begin{equation*}
\begin{cases}
y=-\partial_1 h(x,x'),\\
y'=\partial_2 h(x,x'),
\end{cases}
\end{equation*}
where $F(x,y)=(x',y')$. The lift $F$ gives rise to a dynamical
system whose orbits are given by the images of points of $\R^2$
under the successive iterates of $F$. The orbit of the point
$(x_0,y_0)$ is the bi-infinite sequence
\[\{...,(x_{-k},y_{-k}),...,(x_{-1},y_{-1}),(x_0,y_0),(x_1,y_1),...,(x_k,y_k),...\},\]
where $(x_k,y_k)=F(x_{k-1},y_{k-1})$. The sequence
\[(...,x_{-k},...,x_{-1},x_0,x_1,...,x_k,...)\] denoted by $(x_i)_{i\in\Z}$ is called
stationary configuration which stratifies the identity
\[\partial_1 h(x_i,x_{i+1})+\partial_2 h(x_{i-1},x_i)=0,\ \text{for\ every\ }i\in\Z.\]
Given a sequence of points $(z_i,...,z_j)$, we can associate its
action
\[h(z_i,...,z_j)=\sum_{i\leq s<j}h(z_s,z_{s+1}).\] A configuration $(x_i)_{i\in\Z}$
is called minimal if for any $i<j\in \Z$, the segment of
$(x_i,...,x_j)$ minimizes $h(z_i,...,z_j)$ among all segments
$(z_i,...,z_j)$ of the configuration  satisfying $z_i=x_i$ and
$z_j=x_j$. It is easy to see that every minimal configuration is a
stationary configuration. By [B], minimal configurations satisfy a
group of remarkable properties as follows:
\begin{itemize}
\item Two distinct minimal configurations cross at most once, which
is so called Aubry's crossing lemma.
\item For every minimal configuration $\bold{x}=(x_i)_{i\in\Z}$, the limit
\[\rho(\bold{x})=\lim_{n\rightarrow\infty}\frac{x_{i+n}-x_i}{n}\]
exists and doesn't depend on $i\in\Z$. $\rho(\bold{x})$ is called
the rotation number of $\bold{x}$.
\item For every $\omega\in \R$, there exists a minimal configuration
with rotation number $\omega$. Following the notations of [B], the
set of all minimal configurations with rotation number $\omega$ is
denoted by $M_\omega^h$, which can be endowed with the topology
induced from the product topology on $\R^\Z$. If
$\bold{x}=(x_i)_{i\in\Z}$ is a minimal configuration, considering
the projection $pr:\ M_\omega^h\rightarrow\R$ defined by
$pr(\bold{x})=x_0$, we set $\mathcal {A}_\omega^h=pr(M_\omega^h)$.
\item If $\omega\in\Q$, say $\omega=p/q$ (in lowest terms), then it is convenient to define the rotation symbol to detect the structure of
$M_{p/q}^h$. If $\bold{x}$ is a minimal configuration with rotation
number $p/q $, then the rotation symbol $\sigma(\bold{x})$ of
$\bold{x}$ is defined as follows
\begin{equation*}
\sigma(\bold{x})=\left\{\begin{array}{ll}
\hspace{-0.4em}p/q+,&\text{if}\ x_{i+q}>x_i+p\ \text{for\ all\ }i,\\
\hspace{-0.4em}p/q,&\text{if}\ x_{i+q}=x_i+p\ \text{for\ all\ }i,\\
\hspace{-0.4em}p/q-,&\text{if}\ x_{i+q}<x_i+p\ \text{for\ all\ }i.\\
\end{array}\right.
\end{equation*}
 Moreover, we set
\begin{align*}
&M_{{p/q}^+}^h=\{\bold{x} \text{\ a is minimal configuration with
rotation symbol}\  p/q \text{\ or\ } p/q+\},\\
&M_{{p/q}^-}^h=\{\bold{x} \text{\ a is minimal configuration with
rotation symbol}\  p/q \text{\ or\ } p/q-\},
\end{align*}
then both $M_{{p/q}^+}^h$ and $M_{{p/q}^+}^h$ are totally ordered.
Namely, every two configurations in each of them do not cross. We
denote $pr(M_{{p/q}^+}^h)$ and $pr(M_{{p/q}^-}^h)$ by $\mathcal
{A}_{{p/q}^+}^h$ and $\mathcal {A}_{{p/q}^-}^h$ respectively.
\item If $\omega\in\R\backslash\Q$ and $\bold{x}$ is a minimal
configuration with rotation number $\omega$, then
$\sigma(\bold{x})=\omega$ and $M_\omega^h$ is totally ordered.
\item $\mathcal {A}_\omega^h$ is a closed subset of $\R$ for every rotation symbol
$\omega$.
\end{itemize}
\subsection{Peierls's barrier}
In [M3], Mather introduced the notion of Peierls's barrier and gave
a criterion of existence of invariant circle. Namely, the exact
area-preserving monotone twist map generated by $h$ admits an
invariant circle with rotation number $\omega$ if and only if the
Peierls's barrier $P_\omega^h(\xi)$ vanishes identically for all
$\xi\in\R$. The Peierls's barrier is defined as follows:
\begin{itemize}
\item If $\xi\in \mathcal {A}_\omega^h$, we set $P_\omega^h(\xi)$=0.
\item If $\xi \not\in \mathcal {A}_\omega^h$, since $\mathcal {A}_\omega^h$ is a closed set in $\R$, then $\xi$ is contained in some
complementary interval $(\xi^-,\xi^+)$ of $\mathcal {A}_\omega^h$ in
$\R$. By the definition of $\mathcal {A}_\omega^h$, there exist
minimal configurations with rotation symbol $\omega$,
$\bold{x^-}=(x_i^-)_{i\in\Z}$ and $\bold{x^+}=(x_i^+)_{i\in\Z}$
satisfying $x_0^-=\xi^-$ and $x_0^+=\xi^+$. For every configuration
$\bold{x}=(x_i)_{i\in\Z}$ satisfying $x_i^-\leq x_i\leq x_i^+$, we
set
\[G_\omega(\bold{x})=\sum_I(h(x_i,x_{i+1})-h(x_i^-,x_{i+1}^-)),\]
where $I=\Z$, if $\omega$ is not a rational number, and $I=\{0,...,
q-1 \}$, if $\omega=p/q$. $P_\omega^h(\xi)$ is defined as the
minimum of $G_\omega(\bold{x})$ over the configurations $\bold{x}\in
\Pi=\prod_{i\in I}[x_i^-,x_i^+]$ satisfying $x_0=\xi$. Namely
\[P_\omega^h(\xi)=\min_{\bold{x}}\{G_\omega(\bold{x})|\bold{x}\in \Pi\ \text{and}\ \ x_0=\xi\}.\]
\end{itemize}
By [M3], $P_\omega^h(\xi)$ is a non-negative periodic function of
the variable $\xi\in\R$ with the modulus of continuity with respect
to $\omega$.

To our example (see Section $3$), the modulus of continuity of
$P_\omega^h(\xi)$ with respect to $\omega$ can be improved
significantly. Loosely speaking, the hyperbolicity of the
perturbation implies the
 exponential approximation from $P_\omega^h(\xi)$ to $P_{0^+}^{h}(\xi)$.
 The details will be provided in Section $5$.

\section{\sc Construction of the generating functions}

Consider a completely integrable system with the generating function
\[h_0(x,x')=\frac{1}{2}(x-x')^2 \quad x,x'\in \R.\] We construct the perturbation consisting of two
parts. The first one is
\begin{equation}\label{31}
u_n(x)=\frac{1}{n^a}(1-\cos(2\pi x) )\quad x\in \R,\end{equation}
where $n\in \N$ and $a$ is a positive constant independent of $n$.
The second one is a non negative function $v_n(x)$ satisfying
\begin{equation}\label{32}
\begin{cases}
 v_n(x+1)=v_n(x),\\
\text{supp}\,v_n\cap [0,1]\subset
[\frac{1}{2}-\frac{1}{n^a},\frac{1}{2}+\frac{1}{n^a}], \\
\max v_n=n^{-s}, \\
{||v_n||}_{C^k}=O(n^{-s'}),
\end{cases}
\end{equation}
where we require $s'>a$. It is enough to take $s=(k+2)a$ for
achieving that. The generating function of the nearly integrable
system is constructed as follow:
\begin{equation}\label{h}
h_n(x,x')=h_0(x,x')+u_n(x')+v_n(x'),
\end{equation}
where $n\in\N$. Moreover, we have the following theorem.

\bthm\label{MR} For $\omega\in\R\backslash\Q$ and $n$ large enough,
the exact area-preserving monotone twist map generated by $h_n$ does
not admit any invariant circles with the rotation number satisfying
\[|\omega|<n^{-\frac{a}{2}-\delta}, \] where $\delta$ is a small positive constant independent of $n$. \ethm
We will prove Theorem \ref{MR} in the following sections. First of
all, based on the theorem, we verify that our example has the
property aforementioned in Section 1.

If $\omega\in \Q$, then the invariant circles with rotation number
$\omega$ could be easily destructed even though the perturbation is
$C^\infty$ close to $0$. Therefore it suffices to consider the
irrational $\omega$. The case with a given irrational rotation
number can be easily reduced to the one with a small enough rotation
number. More precisely,

\blm\label{Herm} Let $h_P$ be a generating function as follow
\[h_P(x,x')=h_0(x,x')+P(x'),\] where $P$ is a periodic
function of periodic $1$. Let $Q(x)=q^{-2}P(qx),q\in \N$, then the
exact area-preserving monotone twist map generated by
$h_Q(x,x')=h_0(x,x')+Q(x')$ admits an invariant circle with rotation
number $\omega \in \R\backslash \Q$ if and only if the exact
area-preserving monotone twist map generated by $h_P$ admits an
invariant circle with rotation number $q\omega-p, p\in \Z$. \elm

We omit the proof and for more details, see [H2]. For the sake of
simplicity of notations, we denote $Q_{q_n}$ by $Q_n$ and the same
to $u_{q_n}, v_{q_n}$ and $h_{q_n}$. Let
\[Q_n(x)={q_n}^{-2}(u_n(q_nx)+v_n(q_nx)),\] where $(q_n)_{n\in \N}$
is a sequence satisfying Dirichlet approximation
\begin{equation}\label{diri}
|q_n\omega-p_n|<\frac{1}{q_n}, \end{equation} where $p_n\in \Z$ and
$q_n\in \N$. Since $\omega\in\R\backslash\Q$, we say
$q_n\rightarrow\infty$ as $n\rightarrow\infty$. Let
$\tilde{h}_n(x,x')=h_0(x,x')+Q_n(x')$, we have

\begin{Corollary}\label{Mcor}
For a given rotation number $\omega\in \R\backslash\Q$ and every
$\eps$, there exists $N$ such that for $n>N$, the exact
area-preserving monotone map generated by $\tilde{h}_n$ admits no
invariant circle with rotation number $\omega$ and
\[||\tilde{h}_n-h_0||_{C^{4-\delta'}}<\eps,\]where $\delta'$ is a small positive constant independent of $n$.
\end{Corollary}

\Proof Based on Theorem \ref{MR} and Dirichlet approximation
(\ref{diri}), it suffices to take
\[\frac{1}{q_n}\leq\frac{1}{{q_n}^{\frac{a}{2}+\delta}},\]
which implies \begin{equation}\label{33}
a\leq2-2\delta.\end{equation} From (\ref{31}), (\ref{32}) and
(\ref{h}), it follows that
\begin{align*}
 ||\tilde{h}_n&(x,x')-h_0(x,x')||_{C^r}\\
 &=||Q_n(x')||_{C^r},\\
&\leq{q_n}^{-2}(||u_n(q_nx')||_{C^r}+||v_n(q_nx')||_{C^r}),\\
&\leq{q_n}^{-2}({q_n}^{-a}(2\pi)^r{q_n}^r+C_1{q_n}^{-s'}{q_n}^r),\\
&\leq C_2{q_n}^{r-a-2},
\end{align*}
where $C_1, C_2$ are positive constants only depending on $r$.

To complete the proof, it is enough to make $r-a-2<0$, which
together with (\ref{33}) implies
\[r<a+2\leq 4-2\delta.\] We set ${\delta}'=2\delta$, then the
proof of Corollary \ref{Mcor} is completed.\End

The following sections are devoted to prove Theorem \ref{MR}. For
simplicity, we don't distinguish the constant $C$ in following
different estimate formulas.

\section{\sc Estimate of lower bound of $P_{0^+}^{h_n}$ }

In this section, we will estimate the lower bound of $P_{0^+}^{h_n}$
at the given point. To achieve that, we need to estimate the
distances of pairwise adjacent elements of the minimal
configuration.
\subsection{A spacing lemma}

\blm\label{lowstep} Let $(x_i)_{i\in \Z}$ be a minimal configuration
of $\bar{h}_n$ with rotation symbol $0^+$, then
\[x_{i+1}-x_i=O(n^{-\frac{a}{2}}),\quad \text{for}\quad x_i\in [\frac{1}{4},\frac{3}{4}],\] where
$\bar{h}_n(x_i,x_{i+1})=h_0(x_i,x_{i+1})+u_n(x_{i+1})$. \elm

\Proof Without loss of generality, we assume $x_i\in [0,1]$ for all
$i\in\Z$. By Aubry's crossing lemma, we have
\[0<...<x_{i-1}<x_i<x_{i+1}<...<1.\]We consider the configuration
$(\xi_i)_{i\in \Z}$ defined by
\begin{equation*}
\xi_j= \left\{\begin{array}{ll}\hspace{-0.4em}x_j,& j<i,\\
\hspace{-0.4em}x_{j+1},& j\geq i.\\
\end{array}\right.
\end{equation*}
Since $(x_i)_{i\in \Z}$ is minimal, we have
\[\sum_{i\in \Z}\bar{h}_n(\xi_i,\xi_{i+1})-\sum_{i\in \Z}\bar{h}_n(x_i,x_{i+1})\geq 0.\]
By the definitions of $\bar{h}_n$ and $(\xi_i)_{i\in\Z}$, we have
\begin{align*}
0&\leq\sum_{i\in \Z}\bar{h}_n(\xi_i,\xi_{i+1})-\sum_{i\in
\Z}\bar{h}_n(x_i,x_{i+1})\\
&=\bar{h}_n(x_{i-1},x_{i+1})-\bar{h}_n(x_{i-1},x_{i})-\bar{h}_n(x_{i},x_{i+1})\\
&=(x_{i+1}-x_i)(x_i-x_{i-1})-u_n(x_i).
\end{align*}
Moreover,\[u_n(x_i)\leq(x_{i+1}-x_i)(x_i-x_{i-1})\leq\frac{1}{4}(x_{i+1}-x_{i-1})^2.\]
Therefore, \[x_{i+1}-x_{i-1}\geq 2\sqrt{u_n(x_i)}.\] For $x_i\in
[\frac{1}{4},\frac{3}{4}]$, $u_n(x_i)\geq n^{-a}$, hence,
\begin{equation}\label{ls} x_{i+1}-x_{i-1}\geq 2n^{-\frac{a}{2}}.
\end{equation}

On the other hand, we consider another configuration $(\eta_i)_{i\in
\Z}$ defined by
\begin{equation*}
\eta_j= \left\{\begin{array}{ll}\hspace{-0.4em}x_{j+1},& j<i,\\
\hspace{-0.4em}\frac{1}{2}(x_j+x_{j+1}),& j=i,\\
\hspace{-0.4em}x_j,& j> i.
\end{array}\right.
\end{equation*}Based on the minimality of $(x_i)_{i\in \Z}$, following the deduction as similar as that above, we have
\[x_{i+1}-x_i\leq 2\sqrt{u_n(\eta_i)}\leq 2\sqrt{2}n^{-\frac{a}{2}}.\]
Since $(x_i)_{i\in\Z}$ is a stationary configuration, we have
\begin{align*}
x_{i+1}-x_i&=-\partial_1\bar{h}_n(x_i,x_{i+1}),\\
&=\partial_2\bar{h}_n(x_{i-1},x_i),\\
&=x_i-x_{i-1}+u_n'(x_i).
\end{align*}
Since $u_n'(x)=\frac{2\pi}{n^a}\sin(2\pi x)$, it follows from
$(\ref{ls})$ that\[x_{i+1}-x_i\geq Cn^{-\frac{a}{2}}.\] Therefore,
we have
\[x_{i+1}-x_i=O(n^{-\frac{a}{2}}),\quad x_i\in \left[\frac{1}{4},\frac{3}{4}\right]
.\] The proof of Lemma \ref{lowstep} is completed.\End
\subsection{The lower bound of $P_{0^+}^{h_n}$ }
By the definition of $v_n$,
$\text{supp}\,v_n\cap[0,1]\subset[\frac{1}{2}-\frac{1}{n^a},\frac{1}{2}+\frac{1}{n^a}]$
and $v_n(x+1)=v_n(x)$. Let $(x_i)_{i\in \Z}$ be the minimal
configuration of
$\bar{h}_n(x_i,x_{i+1})=h_0(x_i,x_{i+1})+u_n(x_{i+1})$ with rotation
symbol $0^+$ satisfying $x_0=\frac{1}{2}-\frac{1}{n^a}$, then
\[(x_i)_{i\in \Z}\cap\text{supp}v_n=\emptyset.\]
Moreover, for all $i\in\Z$, \[v_n(x_i)=0.\]

 Let $(\xi_i)_{i\in \Z}$
be a minimal configuration of $h_n$ defined by (\ref{h}) with
rotation symbol $0^+$ satisfying $\xi_0=\eta$, where $\eta$
satisfies $v_n(\eta)=\max v_n(x)=n^{-s}$, then
\begin{align*}
\sum_{i\in
\Z}(h_n(&\xi_i,\xi_{i+1})-h_n(\xi_i^-,\xi_{i+1}^-))\\
&\geq v_n(\eta)+\sum_{i\in \Z}\bar{h}_n(\xi_i,\xi_{i+1})-\sum_{i\in
\Z}h_n(\xi_i^-,\xi_{i+1}^-),\\
&\geq v_n(\eta)+\sum_{i\in \Z}\bar{h}_n(x_i,x_{i+1})-\sum_{i\in
\Z}h_n(x_i,x_{i+1}),\\
&=v_n(\eta)-\sum_{i\in \Z}v_n(x_{i+1}),\\
&=v_n(\eta).
\end{align*}
Therefore,
\[P_{0^+}^{h_n}(\eta)=\min_{x_0=\eta}\sum_{i\in
\Z}(h_n(x_i,x_{i+1})-h_n(x_i^-,x_{i+1}^-))\geq v_n(\eta)=n^{-s}.\]We
conclude that there exists a point $\xi\in
[\frac{1}{2}-\frac{1}{n^a},\frac{1}{2}+\frac{1}{n^a}]$ such that
\begin{equation}\label{lowbound}
P_{0^+}^{h_n}(\xi)\geq n^{-s}.
\end{equation}

\section{\sc The approximation from $P_{0^+}^{h_n}$ to
$P_{\omega}^{h_n}$ } In this section, we will prove the improvement
of modulus of continuity of Peierls's barrier based on the
hyperbolicity of $h_n$. Namely

\blm\label{apprx} For every irrational rotation symbol $\omega$
satisfying $0<\omega<n^{-\frac{a}{2}-\delta}$, we have
\[|P_{\omega}^{h_n}(\xi)-P_{0^+}^{h_n}(\xi)|\leq C\exp(-n^{\delta}).\]
where $\xi\in
\left[\frac{1}{2}-\frac{1}{n^a},\frac{1}{2}+\frac{1}{n^a}\right]$
and $\delta$ is a small positive constant independent of  $n$.
 \elm

\subsection{Some counting lemmas}
To prove the lemma, we need to do some preliminary work. First of
all, we count the number of the elements of a minimal configuration
$(x_i)_{i\in\Z}$ with arbitrary rotation symbol $\omega$ in a given
interval. With the method of [F], we can conclude the following
lemma.

\blm\label{count 0} Let $(x_i)_{i\in\Z}$ be a minimal configuration
of $h_n$ with rotation symbol $\omega$,
$J_n=\left[\exp\left(-n^{\frac{\delta}{2}}\right),\frac{1}{2}\right]$
and $\Lambda_n=\{i\in \Z|\,x_i\in J_n\}$, then \[\sharp\Lambda_n\leq
Cn^{\frac{a}{2}+\frac{\delta}{2}},\] where $\sharp\Lambda_n$ denotes
the number of elements in $\Lambda_n$ and  $\delta$ is a small
positive constant independent of $n$. \elm

\Proof Let $x^+=1-\exp\left(-n^{\frac{\delta}{2}}\right),
x^-=\frac{1}{2}$ and
$\sigma=\left(\frac{x^+}{x^-}\right)^{\frac{1}{N}}$, hence,
\[\ln \sigma=\frac{\ln(x^+)-\ln(x^-)}{N}.\] We choose $N\in \N$
such that $1\leq\ln\sigma\leq 2$, then
$N=O\left(n^{\frac{\delta}{2}}\right)$.

We consider the partition of the interval $J_n=[x^-,x^+]$ into the
subintervals $J_n^k=[\sigma^k x^-,\sigma^{k+1}x^-]$ where $0\leq
k<N$. Hence, $J_n=\cup_{k=0}^{N-1}J_n^k$. We set
$S_k=\{i\in\Lambda_n|(x_{i-1},x_{i+1})\subset J_n^k\}$ and
$m_k=\sharp S_k$.

By the similar deduction as the one in Lemma \ref{lowstep}, we have
\[x_{i+1}-x_{i-1}\geq 2\sqrt{u_n(x_i)+v_n(x_i)}\geq Cn^{-\frac{a}{2}}x_i,\quad\text{for}\quad x_i\in \left[0,\frac{1}{2}\right].\]
For simplicity of notation, we write $Cn^{-\frac{a}{2}}$ by
$\alpha_n$.

If  there exists $k$ such that $i\in S_k$ for $(x_i)_{i\in \Z}$,
then $x_{i+1}-x_{i-1}\geq\alpha_n\sigma^kx^-$,
moreover,\[m_k\alpha_n\sigma^kx^-\leq 2\mathcal {L}(J_n^k)=
2(\sigma-1)\sigma^kx^-,\]where $\mathcal {L}(J_n^k)$ denotes the
length of the interval of $J_n^k$. Hence $m_k\leq 2(\sigma
-1)\alpha_n^{-1}$.

On the other hand, if $i\not\in S_k$ for any $k$, then there exists
$l$ satisfying $1\leq l<N$ such that
\[x_{i-1}<\sigma^lx^-<x_{i+1}.\]Hence,
\[\sharp\{i\in\Lambda_n|i\not\in S_k \ \text{for\ any}\ k\}\leq 2N.
\] Therefore,
\begin{align*}
\sharp(\Lambda_n)\leq 2N(\sigma-1)\alpha_n^{-1}+2N.
\end{align*}Since $1\leq\ln\sigma\leq 2$ and $N=O\left(n^{\frac{\delta}{2}}\right)$, then we have
\[\sharp\Lambda_n\leq
Cn^{\frac{a}{2}+\frac{\delta}{2}}.\] The proof of Lemma \ref{count
0} is completed.
\begin{Remark}
Let $(x_i)_{i\in\Z}$ be a minimal configuration of $h_n$ defined by
$(\ref{h})$ with rotation symbol $\omega$, An argument as similar as
the one in Lemma \ref{count 0} implies that
\[\sharp\left\{i\in\Z\bigg|\,x_i\in\left[\exp\left(-n^{\frac{\delta}{2}}\right),1-\exp\left(-n^{\frac{\delta}{2}}\right)\right]\right\}\leq Cn^{\frac{a}{2}+\frac{\delta}{2}}.\]
\end{Remark}

Second, it is easy to count the number of the elements of a minimal
configuration with irrational rotation symbol. More precisely, we
have the following lemma.
 \blm\label{count w}Let
$(x_i)_{i\in \Z}$ be a minimal configuration with rotation number
$\omega\in \R\backslash\Q$. Then for every interval $I_k$ of length
$k$, $k\in \N$,
\[\frac{k}{\omega}-1\leq\sharp\{i\in \Z|x_i\in
I_k\}\leq\frac{k}{\omega}+1.\] \elm \Proof For every minimal
configuration $(x_i)_{i\in\Z}$ with rotation number $\omega$, there
exists an orientation-preserving circle homeomorphism $\phi$ such
that $\rho(\Phi)=\omega$, where $\Phi:\R\rightarrow\R$ denotes the
lift of $\phi$. Since $\omega\in \R\backslash\Q$, thanks to [H1],
$\phi$ has a unique invariant probability measure $\mu$ on $\T$ such
that $\mu[x,\Phi(x)]=\omega$ for every $x\in\R$. In particular,
\[\mu[x_i,x_{i+1}]=\omega,\quad \text{for\ every\ }i\in\Z.\]
From $\mu(I_k)=k$, it follow that
\begin{align*}
&\omega(\sharp\{i\in \Z|x_i\in I_k\}-1)\leq k,\\
&\omega(\sharp\{i\in \Z|x_i\in I_k\}+1)\geq k,
 \end{align*}
which completes the proof of Lemma \ref{count w}.\End

Based on Lemma \ref{count 0} and Lemma \ref{count w}, if
$0<\omega<n^{-\frac{a}{2}-\delta}$ and $\omega$ is irrational, then
\[\sharp\{i\in \Z|x_i\in [0,1]\}\geq \frac{1}{\omega}-1\geq
C_1n^{\frac{a}{2}+\delta}>C_2n^{\frac{a}{2}+\frac{\delta}{2}}.\] So
far, we have proved the following conclusion.

\blm\label{lbspace} Let $(x_i)_{i\in \Z}$ be a minimal configuration
of $h_n$ defined by $(\ref{h})$ with rotation symbol
$0<\omega<n^{-\frac{a}{2}-\delta}$, then there exists $j^-,j^+\in\Z$
such that \begin{align*}&0<x_{j^--1}<x_{j^-}<x_{j^-+1}\leq
\exp(-n^{\frac{\delta}{2}}),\\
& 1-\exp(-n^{\frac{\delta}{2}})\leq
x_{j^+-1}<x_{j^+}<x_{j^++1}<1.\end{align*}\elm Without loss of
generality, we assume that \begin{equation}\label{jj}j^+-j^-\geq
Cn^{\frac{a}{2}+\frac{2\delta}{3}}.\end{equation}

 If $\xi\in
\mathcal {A}_\omega^{h_n}$, then $P_\omega^{h_n}(\xi)=0$. Hence, it
suffices to consider the case with $\xi\not\in \mathcal
{A}_\omega^{h_n}$ for destruction of invariant circles.  Let
$(\xi^-,\xi^+)$ be the complementary interval of $\mathcal
{A}_\omega^{h_n}$ in $\R$ and contains $\xi$. Let
$\pmb{\xi^{\pm}}=(\xi_i^{\pm})_{i\in\Z}$ be the minimal
configurations with rotation symbol $\omega$ satisfying
$\xi_0^{\pm}=\xi^{\pm}$ and  let $(\xi_i)_{i\in \Z}$ be a minimal
configuration of $h_n$ with rotation symbol $\omega$ satisfying
$\xi_0=\xi$ and $\xi_i^-\leq \xi_i\leq \xi_i^+$. By the definition
of Peierls barrier, we have
\[P_\omega^{h_n}(\xi)=\sum_{i\in \Z}(h_n(\xi_i,\xi_{i+1})-h_n(\xi_i^-,\xi_{i+1}^-)).\]
Since $P_\omega^{h_n}(\xi)$ is 1-periodic with respect to $\xi$,
without loss of generality, we assume that $\xi\in [0,1]$. We set
$d(x)=\min\{|x|,|x-1|\}$ and write $\exp(-n^{\frac{\delta}{2}})$ by
$\epsilon(n)$. By Lemma \ref{lbspace}, there exist $i^-,\ i^+$ such
that
\begin{equation}\label{5}
d(\xi_i^-)<\epsilon(n)\quad \text{and}\quad
\xi_{i+1}^--\xi_{i-1}^-\leq\epsilon(n)\quad \text{for}\quad i=i^-,\
i^+. \end{equation}
 Thanks to Aubry's crossing lemma, we have $\xi_i^-\leq
\xi_i\leq \xi_i^+\leq \xi_{i+1}^-$. Hence,
\[\xi_i-\xi_i^-\leq\epsilon(n)\quad \text{for}\quad i=i^-,\ i^+.\]

\subsection{Proof of lemma \ref{apprx}}
In the following, we will prove Lemma \ref{apprx} with the method
similar to the one developed by Mather in [M3]. The proof can be
proceeded in the following two steps.
 \subsubsection{Step 1} We consider the number of the elements in a segment of the configuration as the length of the segment. In the first step,
 we approximate
$P_\omega^{h_n}(\xi)$ for $\xi\in
\left[\frac{1}{2}-\frac{1}{n^a},\frac{1}{2}+\frac{1}{n^a}\right]$ by
the difference of the actions of the segments of length $i^+-i^-+1$.
To achieve that, we define the following configurations
\begin{align*}
x_i=\left\{\begin{array}{ll}\hspace{-0.4em} \xi_i,& i\neq i^-,\ i^+,\\
\hspace{-0.4em}\xi_i^-,& i=i^-,\ i^+,
\end{array}\right.
\quad \text{and}\ \quad y_i=\left\{\begin{array}{ll}
\hspace{-0.4em}\xi_i,& i^-<i<i^+,\\
\hspace{-0.4em}\xi_i^-,& i\leq i^-,\ i\geq i^+.
\end{array}\right.
\end{align*}
It is easy to see that $\xi_0=\xi$ is contained both of
$(x_i)_{i\in\Z}$ and $(y_i)_{i\in\Z}$ up to the rearrangement of the
index $i$ since $\xi\in
\left[\frac{1}{2}-\frac{1}{n^a},\frac{1}{2}+\frac{1}{n^a}\right]$.
Hence, by the minimality of $(\xi_i)_{i\in \Z}$ satisfying
$\xi_0=\xi$, we have
\begin{equation}\label{c}
P_\omega^{h_n}(\xi)\leq\sum_{i\in
\Z}(h_n(y_i,y_{i+1})-h_n(\xi_i^-,\xi_{i+1}^-)). \end{equation} Since
$\omega$ is irrational, then $(x_i)_{i\in \Z}$ is asymptotic to
$(\xi_i^-)_{i\in \Z}$, which together with the minimality of
$(\xi_i^-)_{i\in \Z}$ yields
\begin{equation}\label{b}
\sum_{i\in
\Z}(h_n(y_i,y_{i+1})-h_n(\xi_i^-,\xi_{i+1}^-)\leq\sum_{i\in
\Z}(h_n(x_i,x_{i+1})-h_n(\xi_i^-,\xi_{i+1}^-)). \end{equation} We
set
\[h(x_i,...,x_j)=\sum_{i\leq s<j}h(x_s,x_{s+1}),\] then
\[\sum_{i\in
\Z}(h_n(x_i,x_{i+1})-h_n(\xi_i,\xi_{i+1}))=\sum_{i=i^-,i^+}(h_n(\xi_{i-1},\xi_i^-,\xi_{i+1})-h_n(\xi_{i-1},\xi_i,\xi_{i+1})).\]
 By the construction of $v_n$ and Lemma \ref{lbspace}, we have $v_n(\xi_{i^-}),\ v_n(\xi_{i^-}^-)=0$. It
follows that
\begin{align*}
h_n(&\xi_{i^--1},\xi_{i^-}^-, \xi_{i^-+1})-h_n(\xi_{i^--1},\xi_{i^-},\xi_{i^-+1})\\
&=h_n(\xi_{i^--1},\xi_{i^-}^-)+h_n(\xi_{i^-}^-,\xi_{i^-+1})-h_n(\xi_{i^--1},\xi_{i^-})-h_n(\xi_{i^-},\xi_{i^-+1}),\\
&=(\xi_{i^-}-\xi_{i^-}^-)(\xi_{i^--1}+\xi_{i^-}^-+\xi_{i^-}+\xi_{i^-+1})+u_n(\xi_{i^-})-u_n(\xi_{i^-}^-),\\
&\leq 4(\xi_{i^-}-\xi_{i^-}^-)\epsilon(n)+u_n'(\eta)(\xi_{i^-}-\xi_{i^-}^-),\\
&\leq 4\epsilon(n)^2+\frac{2\pi}{n^a}\sin(2\pi\eta)\epsilon(n),\\
&\leq C\epsilon(n)^2,
\end{align*}
where $\eta\in (\xi_{i^-},\xi_{i^-}^-)$. It is similar to obtain
\[h_n(\xi_{i^+-1},\xi_{i^+}^-,\xi_{i^++1})-h_n(\xi_{i^+-1},\xi_{i^+},\xi_{i^++1})\leq C\epsilon(n)^2.\]
Hence,
\begin{equation}\label{a}
\sum_{i\in \Z}(h_n(x_i,x_{i+1})-h_n(\xi_i,\xi_{i+1}))\leq
C\epsilon(n)^2.\end{equation} Moreover,
\begin{align*}
\sum_{i\in \Z}(h_n&(x_i,x_{i+1})-h_n(\xi_i^-,\xi_{i+1}^-))\\
&=\sum_{i\in
\Z}(h_n(x_i,x_{i+1})-h_n(\xi_i,\xi_{i+1})+h_n(\xi_i,\xi_{i+1})-h_n(\xi_i^-,\xi_{i+1}^-)),\\
&=\sum_{i\in
\Z}(h_n(x_i,x_{i+1})-h_n(\xi_i,\xi_{i+1}))+P_\omega^{h_n}(\xi),\\
&\leq P_\omega^{h_n}(\xi)+C\epsilon(n)^2.
\end{align*}
Therefore, it follows from $(\ref{c})$ and $(\ref{b})$ that
\begin{equation}\label{step1}
P_\omega^{h_n}(\xi)\leq\sum_{i\in
\Z}(h_n(y_i,y_{i+1})-h_n(\xi_i^-,\xi_{i+1}^-))\leq
P_\omega^{h_n}(\xi)+C\epsilon(n)^2,
\end{equation}
where
\[\sum_{i\in
\Z}(h_n(y_i,y_{i+1})-h_n(\xi_i^-,\xi_{i+1}^-))=h_n(y_{i^-},...,y_{i^+})-h_n(\xi_{i^-}^-,...,\xi_{i^+}^-).\]

\subsubsection{Step 2}
It follows from [M4] that the Peierls's barrier $P_{0^+}^{h_n}(\xi)$
could be defined as follows
\begin{equation}\label{t}
P_{0^+}^{h_n}(\xi)=\min_{\eta_0=\xi}\sum_{i\in
\Z}h_n(\eta_i,\eta_{i+1})-\min\sum_{i\in \Z}h_n(z_i,z_{i+1}),
\end{equation}
where $(\eta_i)_{i\in\Z}$ and  $(z_i)_{i\in\Z}$ are monotone
increasing configurations limiting on $0,\ 1$. We set
\begin{equation*}
\begin{cases}
K(\xi)=\min_{\eta_0=\xi}\sum_{i\in \Z}h_n(\eta_i,\eta_{i+1}),\\
K=\min\sum_{i\in \Z}h_n(z_i,z_{i+1}).
\end{cases}
\end{equation*}

First of all, it is easy to see that $K(\xi)$ and $K$ are bounded.
Second, $P_{0^+}^{h_n}(\xi)=0$ for $\xi=0$ or $1$. Hence, we only
need to consider the case with $\xi\in (0,1)$. Following the ideas
of [M6], let $\pmb{\xi^-}$ and $\pmb{\xi^+}$ be minimal
configurations of rotation symbol $0^+$ and let $(\xi_0^-,\xi_0^+)$
be the complementary interval of $\mathcal {A}_{0^+}^{h_n}$ and
contains $\xi$. Based on the definition
\[P_{0^+}^{h_n}(\xi)=\min_{x_0=\xi}\{G_{0^+}(\bold{x})|\xi_i^-\leq x_i\leq\xi_i^+\},\]
where \[G_{0^+}(\bold{x})=\sum_{i\in
\Z}(h_n(x_i,x_{i+1})-h_n(\xi_i^-,\xi_{i+1}^-))=-K+\sum_{i\in
\Z}h_n(x_i,x_{i+1}),\] the proof of $(\ref{t})$ will be completed
when we verify that the configuration $(x_i)_{i\in \Z}$ achieving
the minimum in the definition of $K(\xi)$ satisfies $\xi_i^-\leq
x_i\leq\xi_i^+$. It can be easily obtained by Aubry's crossing
lemma. In fact, since $(\xi_i^-)_{i\in\Z}$ and $(x_i)_{i\in\Z}$ are
minimal and both are $\alpha$-asymptotic to $0$ as well as
$\omega$-asymptotic to $1$, by Aubry's crossing lemma,
$(\xi_i^-)_{i\in\Z}$ and $(x_i)_{i\in\Z}$ do not cross. Similarly
$(\xi_i^+)_{i\in\Z}$ and $(x_i)_{i\in\Z}$ do not cross. It follows
from $x_0\in (\xi_0^-,\xi_0^+)$ that $(x_i)_{i\in \Z}$ achieving the
minimum in the definition of $K(\xi)$ satisfies $\xi_i^-\leq
x_i\leq\xi_i^+$.

In the following, we will compare $K$, $K(\xi)$ with
$h_n(\xi_{i^-}^-,...,\xi_{i^+}^-)$, $h_n(y_{i^-},...,y_{i^+})$
respectively.

First, we consider $K$ and $h_n(\xi_{i^-}^-,...,\xi_{i^+}^-)$. Let
$(z_i)_{i\in\Z}$ be a monotone increasing configuration limiting on
$0,\ 1$ such that $K=\sum_{i\in \Z}h_n(z_i,z_{i+1})$. By Lemma
\ref{count 0},
\[\sharp\{i\in
\Z|(z_i)_{i\in \Z}\cap[\epsilon(n),1-\epsilon(n)]\}\leq
Cn^{\frac{a}{2}+\frac{\delta}{2}}.\] On the other hand, since
$(z_i)_{i\in\Z}$ has the rotation number $0^+$, then from
(\ref{jj}), it follows that up to the rearrangement of the index
$i$, there exists a subset of length $i^+-i^-$ of $(z_i)_{i\in\Z}$,
denoted by $\{z_{i^-},z_{i^-+1},\ldots,z_{i^+-1},z_{i^+}\}$ such
that
\[z_{i^-+1}\leq\epsilon(n),\quad z_{i^+-1}\geq\epsilon(n).\]

We consider the configuration $(\bar{x}_i)_{i\in\Z}$ defined by
\begin{equation*}
\left\{\begin{array}{ll}
\hspace{-0.4em}\bar{x}_i=\xi_i^-,&i^-<i<i^+,\\
\hspace{-0.4em}\bar{x}_i=0,& i\leq i^-,\\
\hspace{-0.4em}\bar{x}_i=1,& i\geq i^+.
\end{array}\right.
\end{equation*}
By the definition of $h_n$, $h_n(\bar{x}_i,\bar{x}_{i+1})=0$ for
$i<i^-$ or $i \geq i^+$, then
\[\sum_{i\in\Z}h_n(\bar{x}_i,\bar{x}_{i+1})=h_n(\bar{x}_{i^-},...,\bar{x}_{i^+}).\]
Moreover, by the minimality of $(z_i)_{i\in\Z}$, we have
\begin{equation}\label{aa}
K\leq
\sum_{i\in\Z}h_n(\bar{x}_i,\bar{x}_{i+1})=h_n(\bar{x}_{i^-},...,\bar{x}_{i^+}).
\end{equation}
By the construction of $h_n$, we have
$v_n(\bar{x}_{i^-+1})=v_n(\xi_{i^-+1}^-)=0$. Hence,
\begin{equation}\label{55}
\begin{split}
h_n(&\bar{x}_{i^-},\bar{x}_{i^-+1})-h_n(\xi_{i^-}^-,\xi_{i^-+1}^-)\\
&=\frac{1}{2}(\bar{x}_{i^-}-\bar{x}_{i^-+1})^2+u_n(\bar{x}_{i^-+1})
-\frac{1}{2}(\xi_{i^-}^--\xi_{i^-+1}^-)^2-u_n(\xi_{i^-+1}^-),\\
&=\frac{1}{2}(\xi_{i^-+1}^-)^2-\frac{1}{2}(\xi_{i^-+1}^--\xi_{i^-}^-)^2,\\
&=\frac{1}{2}\xi_{i^-}^-(2\xi_{i^-+1}^--\xi_{i^-}^-),\\
&\ \leq C\epsilon(n)^2.
\end{split}
\end{equation}
It is similar to obtain
\begin{equation}\label{3}
h_n(\bar{x}_{i^+-1},\bar{x}_{i^+})
-h_n(\xi_{i^+-1}^-,\xi_{i^+}^-)\leq C\epsilon(n)^2. \end{equation}
 Since
\begin{align*}
h_n(&\bar{x}_{i^-},...,\bar{x}_{i^+})-h_n(\xi_{i^-}^-,...,\xi_{i^+}^-)\\
&=h_n(\bar{x}_{i^-},\bar{x}_{i^-+1})-h_n(\xi_{i^-}^-,\xi_{i^-+1}^-)+h_n(\bar{x}_{i^+-1},\bar{x}_{i^+})
-h_n(\xi_{i^+-1}^-,\xi_{i^+}^-),
\end{align*}then
\begin{equation}\label{aaa}
h_n(\bar{x}_{i^-},...,\bar{x}_{i^+})-h_n(\xi_{i^-}^-,...,\xi_{i^+}^-)\leq
C\epsilon(n)^2.
\end{equation}
 From $(\ref{aa})$ and $(\ref{aaa})$ we have
\begin{equation}\label{1}
 K\leq h_n(\xi_{i^-}^-,...,\xi_{i^+}^-)+C\epsilon(n)^2.
\end{equation}

To obtain the reverse inequality of $(\ref{1})$, we consider the
configuration as follows
\begin{equation*}
\left\{\begin{array}{ll}
\hspace{-0.4em}\tilde{x}_i=z_i,&i^-<i<i^+,\\
\hspace{-0.4em}\tilde{x}_i=0,&i\leq i^-,\\
\hspace{-0.4em}\tilde{x}_i=1,&i\geq i^+.
\end{array}\right.
\end{equation*}

From the definition of $h_n$, it follows that $v_n(z_{i^-+1})=0$ and
$h_n(z_i,z_{i+1})\geq 0$ for all $i\in \Z$. Moreover, we have
\begin{align*}
h_n(&\tilde{x}_{i^-},...,\tilde{x}_{i^+})-K\\
&=h_n(\tilde{x}_{i^-},\tilde{x}_{i^-+1})+h_n(\tilde{x}_{i^+-1},\tilde{x}_{i^+})-\sum_{i<i^-,i\geq
i^+}h_n(z_i,z_{i+1}),\\
&\leq\frac{1}{2}(z_{i^-+1})^2+u_n(z_{i^-+1})+\frac{1}{2}(z_{i^+-1}-1)^2,\\
&\leq u_n'(\eta)z_{i^-+1}+C_1\epsilon(n)^2,\\
&\leq 2\pi n^{-a}\sin(2\pi\eta)z_{i^-+1}+C_1\epsilon(n)^2,\\
&\leq C_2n^{-a}(z_{i^-+1})^2+C_1\epsilon(n)^2,\\
&\leq C\epsilon(n)^2.
\end{align*}
where $\eta\in (0,z_{i^-+1})$. Namely
\begin{equation}\label{bb}
h_n(\tilde{x}_{i^-},...,\tilde{x}_{i^+})\leq K+C\epsilon(n)^2.
\end{equation}
Furthermore, we consider the finite segment of the configuration
defined by
\begin{equation*}
\left\{\begin{array}{ll}
\hspace{-0.4em}\eta_i=\tilde{x}_i,&i^-<i<i^+,\\
\hspace{-0.4em}\eta_i=\xi_i^-,&i=i^-,\\
\hspace{-0.4em}\eta_i=\xi_i^-,& i=i^+.
\end{array}\right.
\end{equation*}
Then, the minimality of $(\xi_i^-)_{i\in\Z}$ implies
$h_n(\xi_{i^-}^-,...,\xi_{i^+}^-)\leq
h_n(\eta_{i^-},...,\eta_{i^+})$. Hence, by $(\ref{bb})$, we have
\begin{equation}\label{cc}
h_n(\xi_{i^-}^-,...,\xi_{i^+}^-)\leq
K+C\epsilon(n)^2+h_n(\eta_{i^-},...,\eta_{i^+})-h_n(\tilde{x}_{i^-},...,\tilde{x}_{i^+}),
\end{equation}
where
\begin{align*}
 h_n(&\eta_{i^-},...,\eta_{i^+})-h_n(\tilde{x}_{i^-},...,\tilde{x}_{i^+})\\
 &=h_n(\eta_{i^-},\eta_{i^-+1})-h_n(\tilde{x}_{i^-},\tilde{x}_{i^-+1})+h_n(\eta_{i^+-1},\eta_{i^+})-h_n(\tilde{x}_{i^+-1},\tilde{x}_{i^+}).
\end{align*}
By the deduction as similar as $(\ref{55})$, we have
\begin{align*}
 &h_n(\eta_{i^-},\eta_{i^-+1})-h_n(\tilde{x}_{i^-},\tilde{x}_{i^-+1})\leq C\epsilon(n)^2,\\
 &h_n(\eta_{i^+-1},\eta_{i^+})-h_n(\tilde{x}_{i^+-1},\tilde{x}_{i^+})\leq C\epsilon(n)^2.
\end{align*}
Moreover, \begin{equation}\label{dd}
 h_n(\eta_{i^-},...,\eta_{i^+})-h_n(\tilde{x}_{i^-},...,\tilde{x}_{i^+})\leq
 C\epsilon(n)^2.
\end{equation}
 Hence, from $(\ref{cc})$ and $(\ref{dd})$, it follows that
\begin{equation}\label{2}
h_n(\xi_{i^-}^-,...,\xi_{i^+}^-)\leq K+C\epsilon(n)^2,
\end{equation}
which together with $(\ref{1})$ and $(\ref{2})$ implies
\begin{equation}\label{13}
|h_n(\xi_{i^-}^-,...,\xi_{i^+}^-)-K|\leq C\epsilon(n)^2.
\end{equation}

Next, we compare $h_n(y_{i^-},...,y_{i^+})$ with $K(\xi)$.

Since $(\xi_i)_{i\in\Z}$ is minimal among all configurations with
rotation symbol $\omega$ satisfying $\xi_0=\xi$. By $(\ref{5})$ and
Aubry's crossing lemma, we have \[d(\xi_i)\leq \epsilon(n),\quad
\text{for}\ i=i^-,i^+,\] where $d(\xi_i)=\min\{|\xi_i|,|\xi_i-1|\}$.
By an argument as similar as the one in the comparison between $K$
and $h_n(\xi_{i^-}^-,...,\xi_{i^+}^-)$, we have
\begin{equation}\label{12}
|h_n(\xi_{i^-},...,\xi_{i^+})-K(\xi)|\leq C\epsilon(n)^2.
\end{equation}
By the construction of $(y_i)_{i\in\Z}$, namely
\begin{equation*}
y_i=\left\{\begin{array}{ll}
\hspace{-0.4em}\xi_i,& i^-<i<i^+,\\
\hspace{-0.4em}\xi_i^-,& i\leq i^-,\ i\geq i^+,
\end{array}\right.
\end{equation*}
we have
\begin{align*}
h_n(&y_{i^-},...,y_{i^+})-h_n(\xi_{i^-},...,\xi_{i^+})\\
&=h_n(\xi_{i^-}^-,\xi_{i^-+1})-h_n(\xi_{i^-},\xi_{i^-+1})+h_n(\xi_{i^+-1},\xi_{i^+}^-)-h_n(\xi_{i^+-1},\xi_{i^+}).
\end{align*}
By the deduction as similar as $(\ref{aaa})$, we have
\begin{equation}\label{11}
|h_n(y_{i^-},...,y_{i^+})-h_n(\xi_{i^-},...,\xi_{i^+})|\leq
C\epsilon(n)^2.
\end{equation}

Finally, from $(\ref{step1})$, $(\ref{13})$, $(\ref{12})$ and
$(\ref{11})$ we obtain
\begin{align*}
|P_{\omega}^{h_n}(\xi)-P_{0^+}^{h_n}(\xi)|&\leq|h_n(y_{i^-},...,y_{i^+})-h_n(\xi_{i^-}^-,...,\xi_{i^+}^-)+K-K(\xi)|+C_1\epsilon(n)^2,\\
&\leq
|h_n(\xi_{i^-},...,\xi_{i^+})-K(\xi)|+|h_n(\xi_{i^-}^-,...,\xi_{i^+}^-)-K|\\
&\ +|h_n(y_{i^-},...,y_{i^+})-h_n(\xi_{i^-},...,\xi_{i^+})|+C_1\epsilon(n)^2,\\
&\leq C\epsilon(n)^2,\\
&=C\exp(-n^\delta),
\end{align*}
which completes the proof of Lemma \ref{apprx}.\End

\section{\sc Proof of theorem \ref{MR}}
Based on the preparation above, it is easy to prove Theorem
\ref{MR}. We assume that there exists an invariant circle with
rotation number $0<\omega<n^{-\frac{a}{2}-\delta}$ for $h_n$, then
$P_\omega^{h_n}(\xi)\equiv 0$ for every $\xi\in \R$. By Lemma
\ref{apprx}, we have
\begin{equation}\label{pp}
|P_{0^+}^{h_n}(\xi)|\leq C\exp(-n^{\delta}),\quad\text{for}\quad
\xi\in
\left[\frac{1}{2}-\frac{1}{n^a},\frac{1}{2}+\frac{1}{n^a}\right].
\end{equation}

On the other hand, $(\ref{lowbound})$ implies that there exists a
point $\tilde{\xi}\in
\left[\frac{1}{2}-\frac{1}{n^a},\frac{1}{2}+\frac{1}{n^a}\right]$
such that
\[P_{0^+}^{h_n}(\tilde{\xi})\geq n^{-s}.\]Hence, we have
\[n^{-s}\leq C\exp(-n^{\delta}).\] It is an obvious
contradiction for $n$ large enough. Therefore, there exists no
invariant circle with rotation number
$0<\omega<n^{-\frac{a}{2}-\delta}$.

For $-n^{-\frac{a}{2}-\delta}<\omega<0$, by comparing
$P_{\omega}^{h_n}(\xi)$ with $P_{0^-}^{h_n}(\xi)$, the proof is
similar. We omit the details. Therefore, the proof of Theorem
\ref{MR} is completed. \End

\noindent\textbf{Acknowledgement} The author would like to thank
Prof. C-Q.Cheng, W.Cheng for many helpful discussions, and X-J. Cui
for his careful proofreading. This work is under the support of the
National Basic Research Programme of China (973 Programme,
2007CB814800) and Basic Research Programme of Jiangsu Province,
China (BK2008013).

{\sc Department of Mathematics, Nanjing University, Nanjing 210093,
China.}

{\it E-mail address:} \texttt{linwang.math@gmail.com}

\end{document}